УДК 621. 372. 061 9 (045)

**А. В. Васильев,** канд. техн. наук
Отделение гибридных моделирующих и управляющих систем в энергетике ИПМЭ им. Г. Е. Пухова НАН Украины, г. Киев

# МЕТОД ЧИСЛЕННО-АНАЛИТИЧЕСКОГО РЕШЕНИЯ ЛИНЕЙНЫХ ДИФФЕРЕНЦИАЛЬНЫХ УРАВНЕНИЙ С ПРОИЗВОДНЫМИ НЕЦЕЛОГО ПОРЯДКА ПО КАПУТО И С ПЕРЕМЕННЫМИ КОЭФФИЦИЕНТАМИ

Рассмотрено применение аппроксимационно-операционного подхода к решению линейных дифференциальных уравнений дробного порядка с переменными коэффициентами. Показано, что метод может применяться также к решению дифференциальных уравнений как дробного, так и смешанного порядков. Вычислительные эксперименты выполнены в программной среде системы Mathematica®.

**Ключевые слова:** *полиномиальная аппроксимация, дробное исчисление, дифференциальное уравнение нецелого порядка, производная по Капуто, математическое моделирование, динамическая система.*

**Введение**. Исследование динамических систем с иррациональными передаточными функциями и функционирующих во фрактальных средах приводит к необходимости исследовать и решать интегро-дифференциальные уравнения с производными и интегралами нецелых порядков [8]. Получение аналитических решений таких уравнений в большинстве случаев сталкивается с большими трудностями. При численном решении дифференциальных уравнений дробного порядка чаще всего используют формулу Грюнвальда-Летникова [7], которая является обобщением методов конечных разностей и конечных сумм классического математического анализа. В случае использования операционных методов обычно стремятся преобразовать исходное дифференциальное уравнение к интегральному путем интегрирования с порядком интегрального оператора равным порядку старшей производной в заданном уравнении. При этом возникает проблема задания начальных условий для дробных производных различных порядков, численные величины которых при исследовании реальных динамических систем не всегда возможно получить. Избежать необходимости задания начальных условий для дробных производных оказывается возможным благодаря введению определения дробной производной по Капуто [5, 8], для которой достаточно иметь начальные условия для функции и целочисленных производных искомого решения. В данной работе рассмотрен аппрокси-мационный метод решения линейных дифференциальных уравнений дробного порядка с переменными коэффициентами и производными по Капуто. Работа построена следующим образом. В ретроспективном плане приведены известные понятия и формулы метода полиномиальной аппроксимации как операционного, а также основные определения производных и интегралов дробного (нецелого) порядка, используемые в статье. Затем рассмотрено применение метода к получению приближенных решений задачи Коши для линейных дифференциальных уравнений с постоянными и переменными коэффициентами, как целого, так и дробного порядков. Изложение сопровождается иллюстративными примерами, выполненными в среде системы *Mathematica*®.

**Полиномиальная аппроксимация как операционный метод**

Как известно [1, 3], сигнал $x(t)$, заданный на интервале измерения аргумента $0 \leq t < T$, может быть аппроксимирован обобщенным полиномом по некоторой системе базисных (образующих) функций $\vec{S}(t) = [s_1(t), s_2(t), \cdots, s_m(t)]^*$ следующим выражением:

$$x_a(t) = \sum_{i=1}^{m} X_i \cdot s_i(t) = \vec{X}^* \cdot \vec{S}(t), \quad (1)$$

где: $\vec{X} = [X_1, X_2, \cdots, X_m]^*$ — вектор коэффициентов полинома (1), называемый аппроксимирующим полиномиальным спектром сигнала, символом * здесь и далее обозначены операции транспонирования векторов и матриц.

Оптимальная аппроксимация, соответствующая минимуму интеграла квадрата функции ошибки на интервале аппроксимации $J = \int_0^T (x(t) - x_a(t))^2 dt \Rightarrow \min(\vec{X})$, достигается, если аппроксимирующий полиномиальный спектр сигнала находится из следующих выражений:

$$\vec{X} = \mathbf{W}^{-1} \cdot \vec{Q}, \quad (2)$$

$$\mathbf{W} = \begin{bmatrix} w_{11} & w_{12} & \cdots & w_{1m} \\ w_{21} & w_{22} & \cdots & w_{2m} \\ \vdots & \vdots & \ddots & \vdots \\ w_{m1} & w_{m2} & \cdots & w_{mm} \end{bmatrix}, \quad (3)$$

$$w_{ij} = \int_0^T s_i(t) \cdot s_j(t) dt, \quad (4)$$

$$\bar{\mathbf{Q}} = [q_1, q_2, \cdots, q_m]^*, \tag{5}$$

$$q_i = \int_0^T x(t) \cdot s_i(t) dt. \tag{6}$$

Предполагается, что система базисных функций обладает функциональной полнотой, определена на том же интервале изменения аргумента, что и сам сигнал, и интегралы (4) и (6) существуют и могут быть определены численно.

Выражения (1)—(6) трактуются как прямое и обратное операционное преобразования [1], операционный характер которых становится очевидным, если их представить в следующем виде:

$$\bar{\mathbf{X}} = \left( \int_0^T \bar{\mathbf{S}}(t) \cdot \bar{\mathbf{S}}(t)^* dt \right)^{-1} \cdot \left( \int_0^T x(t) \cdot \bar{\mathbf{S}}(t) dt \right), \tag{7}$$

$$x_a(t) = \bar{\mathbf{X}}^* \cdot \bar{\mathbf{S}}(t). \tag{8}$$

Первое из этих выражений (7) сопоставляет сигналу $x(t)$, как функции времени, его операционное изображение в виде вектора коэффициентов аппроксимирующего полинома $\bar{\mathbf{X}}$ и является, по существу, прямым операционным преобразованием. Выражение (8) выполняет реконструкцию сигнала по его изображению в виде аппроксимации и является обратным операционным преобразованием. Необходимо заметить, что операции интегрирования в (7) выполняются над каждым элементом матрицы или вектора. Каждая система базисных функций формирует свой вид операционных преобразований. Каждой математической операции над сигналами соответствуют эквивалентные математические операции над изображениями, которые формируют своеобразную алгебру. Приведем с целью полноты изложения только два положения такой алгебры.

Первое положение касается линейной комбинации двух сигналов:

$$z(t) = a \cdot x(t) \pm b \cdot y(t) \Leftrightarrow \bar{\mathbf{Z}} = a \cdot \bar{\mathbf{X}} \pm b \cdot \bar{\mathbf{Y}}, \tag{9}$$

которое означает, что линейной комбинации сигналов соответствует в области изображений такая же линейная комбинация векторов аппроксимирующих спектров этих сигналов.

Второе положение определяет правило нахождения изображения интеграла сигнала по изображению подинтегральной функции:

$$y(t) = \int_0^t x(\tau) d\tau \Leftrightarrow \bar{\mathbf{Y}} = \mathbf{P}_s^1 \cdot \bar{\mathbf{X}}. \tag{10}$$

В (9), (10) символом $\Leftrightarrow$ обозначено взаимное соответствие выражений в пространстве оригиналов и пространстве изображений. Операции интегрирования сигнала с переменным верхним пределом (10) соответствует умножение вектора аппроксимирующего спектра подинтегральной функции на так называемую операционную матрицу интегрирования $\mathbf{P}_s^1$, элементы которой зависят только от системы базисных функций. Нижний индекс указывает на вид базисной системы, а верхний индекс — порядок интегрального оператора. Систематическое применение операционных преобразований к математической модели динамической системы в форме интегро-дифференциальных уравнений приводит к математической модели той же динамической системы в форме алгебраических уравнений, что позволяет исследовать динамическую систему как алгебраический объект (безинерционную систему). Проиллюстрируем порядок применения операционного подхода к простейшему дифференциальному уравнению первого порядка:

$$\begin{cases} \dfrac{dy(t)}{dt} + a \cdot y(t) = f(t) \\ y(0) = y_0 \end{cases} \qquad (11)$$

Интегрируя (11), получим следующее интегральное уравнение:

$$y(t) + a \cdot \int_0^t y(\tau)d\tau = \int_0^t f(\tau)d\tau + y_0. \qquad (12)$$

Переход в область изображений приводит к следующему векторно-матричному уравнению:

$$\bar{\mathbf{Y}} + a \cdot \mathbf{P}_s^1 \cdot \bar{\mathbf{Y}} = \mathbf{P}_s^1 \cdot \bar{\mathbf{F}} + y_0 \cdot \bar{\mathbf{1}}. \qquad (13)$$

Решение (13) имеет вид:

$$\bar{\mathbf{Y}} = \left(\mathbf{E} + a \cdot \mathbf{P}_s^1\right)^{-1} \cdot (\mathbf{P}_s^1 \cdot \bar{\mathbf{F}} + y_0 \cdot \bar{\mathbf{1}}), \qquad (14)$$

где $\mathbf{E}$ — единичная матрица $m$-го размера.

Переход в область оригиналов приводит к следующей аппроксимации решения:

$$y_a(t) = \bar{\mathbf{Y}}^* \cdot \bar{\mathbf{S}}(t). \qquad (15)$$

Полученное решение является численно-аналитическим, так как в качестве базисной системы функций могут быть использованы различные системы (степенные, экспоненциальные, тригонометрические, ортогональные полиномы и т.д.). Числовой результат формируется при задании конкретной базисной системы.

При аппроксимации сигналов в моделировании и автоматическом управлении широкое распространение получила базисная сис-тема функций, называемых блочно-импульсными [6]. Эта система представляет собой совокупность прямоугольных импульсов единичной амплитуды, сдвинутых относительно друг друга на величину продолжительности импульса и заполняющих интервал аппроксимации. Определение такой базисной системы имеет вид:

$$s_i(t) = v_i(t) = \begin{cases} 1, & \text{если } (i-1)h \leq t < ih \\ 0, & \text{если } t < (i-1)h \text{ или } t \geq ih \end{cases}, \quad i := 1, 2, \cdots m, \quad h = \frac{T}{m}. \quad (16)$$

В такой базисной системе задача аппроксимации сигнала существенно упрощается, и блочно-импульсный спектр сигнала определяется выражением:

$$X_i = \frac{1}{h} \int_{(i-1)h}^{ih} x(t)dt, \qquad (17)$$

а операционная матрица интегрирования имеет вид:

$$\mathbf{P}_v^1 = \frac{h}{2} \cdot \underbrace{\begin{bmatrix} 1 & 0 & \cdots & 0 \\ 2 & 1 & \cdots & 0 \\ \vdots & \vdots & \ddots & \vdots \\ 2 & 2 & \cdots & 1 \end{bmatrix}}_{m \text{ столбцов}}. \qquad (18)$$

Аппроксимация сигналов при использовании блочно-импульсной системы базисных функций является кусочно-ступенчатой, а элементы аппроксимирующего блочно-импульсного спектра сигнала имеют физический смысл среднего значения сигнала на интервале $h$ изменения аргумента. Сравнительно невысокая точность аппроксимации является существенным недостатком блочно-импульсной системы базисных функций. Однако, простота решения задачи аппроксимации, особенно при выполнении нелинейных операций над сигналами, делают эту систему весьма привлекательной и распространенной. В работах [2, 4] предложен интерполяционно-экстраполяционный метод повышения точности аппроксимации, полученной при использовании блочно-импульсных спектров. В данной работе будет использована именно эта система базисных функций.

**Основные определения производных и интегралов нецелых порядков.** Дробное исчисление является обобщением классического математического анализа, имеет уже трехвековую историю и в настоящее время претерпевает второе рождение [5, 7, 8]. Поскольку предметом данной работы является применение операционного подхода к решению дифференциальных уравнений нецелого порядка с переменными коэффициентами, рассмотрим основные определения интегралов и производных дробных порядков, используемые в рабо-те, а также выражение операционных матриц интегрирования дробного порядка в системе блочно-импульсных базисных функций.

### Интеграл Римана-Лиувилля

Наиболее часто применяется определение интеграла дробного порядка, связанное с интегральной формулой Коши и получившее название интеграла Римана-Лиувилля [3, 7, 8]:

$$I^\beta f(t) = \frac{1}{\Gamma(\beta)} \int_0^t (t-\tau)^{\beta-1} f(\tau) d\tau, \qquad (19)$$

где:

$\beta$ — нецелый (дробный) порядок интегрального оператора ($0 < \beta \in R$),

$\Gamma(*)$ — гамма функция Эйлера.

При $\beta = 1$ выражение (19) превращается в обычный интеграл с переменным верхним пределом.

### Производная дробного порядка $\beta$

Производная дробного порядка $\beta$ при значениях порядка, лежащих в пределах $m-1 < \beta \leq m$, $m \in N$, связана с интегралом Римана-Лиувилля следующим выражением:

$$D^\beta f(t) = \frac{d^m}{dt^m} I^{m-\beta} f(t) = \frac{d^m}{dt^m} \frac{1}{\Gamma(m-\beta)} \int_0^t (t-\tau)^{m-\beta-1} f(\tau) d\tau. \quad (20)$$

### Производная дробного порядка по Капуто

В технических приложениях удобным является определение производной дробного порядка, введенное Капуто [3, 8]:

$$^C D^\beta f(t) = \begin{cases} I^{m-\beta} f^{(m)}(t) = \dfrac{1}{\Gamma(m-\beta)} \int_0^t (t-\tau)^{m-\beta-1} \dfrac{d^m f(\tau)}{d\tau^m} d\tau \\ \dfrac{d^m f(t)}{dt^m}, \; \beta = m \end{cases} \qquad (21)$$

При использовании определения производной дробного порядка по Капуто оказывается возможным использование естественных начальных значений функции и целочисленных производных, вместо начальных значений дробных производных различных порядков, которые в большинстве случаев получить затруднительно. При нулевых начальных значениях производные по Капуто совпадают с производными на основе интеграла Римана-Лиувилля.

*Операционная матрица интегрирования дробного порядка*

В системе блочно-импульсных базисных функций операционная матрица интегрирования дробного порядка является нижней треугольной тёплицевой матрицей [3]:

$$\mathbf{P}_V^\beta = \{p_{ij}\}\big|_{i=1, j=1}^{m}. \qquad (22)$$

Ее элементы определяются следующим выражением [3]:

$$p_{ij} = \frac{h^\beta}{\Gamma(\beta+2)} \begin{cases} 0, \text{ если i<j} \\ 1, \text{ если i=j} \\ (i-j+1)^{\beta+1} - 2(i-j)^{\beta+1} + (i-j-1)^{\beta+1}, \text{ если i>j} \end{cases} \qquad (23)$$

Необходимо отметить, что выражение (23) при целых значениях порядка $\beta$ формирует операционную матрицу интегрирования целого порядка, в частности, при $\beta = 1$ получается матрица (18).

### Аппроксимационно-операционный подход к решению дифференциальных уравнений нецелого порядка с переменными коэффициентами и с производными по Капуто

Порядок применения метода рассмотрим на двух иллюстративных примерах.

*Пример 1.* Задана следующая задача Коши для дифференциального уравнения с переменными коэффициентами: Найти аппроксимацию функции $x(t)$, удовлетворяющую дифференциальному уравнению: $^C D^{\beta 1} x(t) + \varphi_1(t) \cdot {^C D^{\beta 2}} x(t) + \varphi_2(t) \cdot x(t) = f(t)$ и начальным условиям: $x(0) = x_{00}$, $\left.\dfrac{dx(t)}{dt}\right|_{t=0} = x_{10}$. Предполагается, что дробные порядки дифференциальных операторов лежат в пределах: $1 < \beta_1 \le 2$, $0 < \beta_2 \le 1$, а переменные коэффициенты удовлетворяют ограничениям: $\varphi_1(t) > 0$, $\varphi_2(t) > 0$.

- На первом этапе решения приводим уравнение к интегральному виду, используя определение дробной производной по Капуто (21):

$$I^{2-\beta 1}\left(\frac{d^2 x(t)}{dt^2}\right) + \varphi_1(t) \cdot I^{1-\beta 2}\left(\frac{dx(t)}{dt}\right) + \varphi_2(t) \cdot x(t) = f(t). \qquad (24)$$

- Вводим для старшей производной целого порядка в (24) новое определение в виде функции $u(t)$:

$$\frac{d^2 x(t)}{dt^2} = u(t). \qquad (25)$$

- Интегрируя (25) и используя начальные условия, формируем выражения для остальных производных целого порядка и самой функции $x(t)$ в терминах введенной функции $u(t)$:

$$\frac{dx(t)}{dt} = x_{10} + \int_0^t u(\tau)d\tau, \qquad (26)$$

$$x(t) = x_{00} + x_{10} \cdot t + \int_0^t\int_0^t u(\tau)d\tau^2. \qquad (27)$$

- Вводим для известных функций и для неизвестных составляющих решения вектора аппроксимирующих спектров:

$$f(t) \Leftrightarrow \bar{F}, \qquad (28)$$

$$\varphi_1(t) \Leftrightarrow \bar{\Phi}_1, \qquad (29)$$

$$\varphi_2(t) \Leftrightarrow \bar{\Phi}_2, \qquad (30)$$

$$u(t) \Leftrightarrow \bar{U}. \qquad (31)$$

- Переходим в область изображений:

$$\begin{aligned}&P^{2-\beta_1}\cdot\bar{U} + D(\bar{\Phi}_1)\cdot P^{1-\beta_2}\cdot P^1\cdot\bar{U} + D(\bar{\Phi}_2)\cdot P^2\cdot\bar{U} = \bar{\Psi},\\ &\bar{\Psi} = \bar{F} - D(\bar{\Phi}_1)\cdot P^{1-\beta_2}\cdot\bar{1}\cdot x_{10} - D(\bar{\Phi}_2)\cdot P^2\cdot\left(x_{00}\cdot\bar{1} + x_{10}\cdot\bar{t}\right)\end{aligned} \qquad (32)$$

В (32) введены диагональные матрицы, элементы которых являются составляющими аппроксимирующих векторов переменных коэффициентов ($\varphi_1(t)$, $\varphi_2(t)$):

$$\mathbf{D}(\bar{\mathbf{Z}}) = \begin{bmatrix} z_1 & 0 & \cdots & 0 \\ 0 & z_2 & \cdots & 0 \\ \vdots & \vdots & \ddots & 0 \\ 0 & 0 & \cdots & z_m \end{bmatrix}, \qquad (33)$$

$$\bar{\mathbf{Z}} = [z_1, z_2, ..., z_m]^*. \qquad (34)$$

- Находим аппроксимирующий спектр функции $u(t)$ путем решения уравнения (32):

$$\bar{U} = \left(P^{2-\beta}D(\bar{\Phi}_1)\cdot P^{2-\beta}D(\bar{\Phi}_2)P^2\right)^{-1}\cdot\Psi. \qquad (35)$$

- Вводим в рассмотрение аппроксимирующие спектры решения уравнения и его первой производной, используя операционные аналоги выражений (26) и (27):

$$\bar{\mathbf{X}}_1 = x_{10}\cdot\bar{1} + \mathbf{P}^1\cdot\bar{U}, \qquad (36)$$

$$\bar{\mathbf{X}} = x_{00}\cdot\bar{1} + x_{10}\cdot\bar{\mathbf{t}} + \mathbf{P}^2\cdot\bar{U}. \qquad (37)$$

- Переходим в область оригиналов и определяем аппроксимации функции решения дифференциального уравнения и ее производных:

$$x_a(t) = \vec{X}^* \cdot \vec{S}(t), \qquad (38)$$

$$\frac{dx_a(t)}{dt} = \vec{X}_1^* \cdot \vec{S}(t), \qquad (39)$$

$$\frac{d^2 x_a(t)}{dt^2} = \vec{U}_1^* \cdot \vec{S}(t). \qquad (40)$$

Приведенный порядок сохраняется и для других дифференциальных уравнений, в том числе целого и смешанного порядков.

Рассмотрим численный пример применения описанного метода.

*Пример 2.* Задана задача Коши для дифференциального уравнения:

$$^C D^{1.5} y(t) + (1+t^2) \cdot {}^C D^{0.3} y(t) + t^2 \cdot x(t) = e^{-t} \qquad (41)$$

при начальных условиях: $y(0) = -5$, $y'(0)=2$ на интервале изменения аргумента $0 \le t < 5$. Необходимо определить аппроксимацию решения уравнения и его производных.

Выберем порядок базисной системы блочно-импульсных функций равным 50. Это соответствует разбиению интервала решения на 50 равных отрезков с длиной $h = T/m = 0.1$. Ниже приводится программа решения задачи предлагаемым методом в среде системы *Mathematica*® [9] и результаты ее работы.

*Программа аппроксимационного решения уравнения (41)*

- *Определение исходных данных:*
```
φ1:=(1+t²); φ1:= t²; f:=e⁻ᵗ
y0=:-5; y1:=2; β1:=1.5; β2:=0.3;
    m:=50; h:=0.1; T:=m*h;
```
- *Определение базисной системы функций:*
```
V:=Table[If[(i-1)*h<t≤i*h,1,0],{i,m}];
```
- *Определение векторов аппроксимирующих спектров известных функций, входящих в уравнение (переменных коэффициентов, правой части, константы 1 и аргумента t):*

$$\Phi 1 := \text{Table}\left[\frac{1}{h} * \int_{(i-1)*h}^{i*h} \varphi 1 \, dt, \{i,m\}\right]$$

$$\Phi 2 := \text{Table}\left[\frac{1}{h} * \int_{(i-1)*h}^{i*h} \varphi 2 \, dt, \{i,m\}\right]$$

$$F := \text{Table}\left[\frac{1}{h} * \int_{(i-1)*h}^{i*h} f \, dt, \{i,m\}\right]$$

$$\text{One} := \text{Table}\left[\frac{1}{h} * \int_{(i-1)*h}^{i*h} 1 \, dt, \{i,m\}\right]$$

$$\text{Tim} := \text{Table}\left[\frac{1}{h} * \int_{(i-1)*h}^{i*h} t \, dt, \{i,m\}\right]$$

- *Определение единичной матрицы и матричных операционных изображений переменных коэффициентов:*

```
E1:=IdentityMatrix[m];
D1:=DiagonalMatrix[Φ1];
D2:=DiagonalMatrix[Φ2];
```

- *Определение операционных матриц интегрирования дробных и необходимых целых порядков:*

```
H[β_,h_,m_]:= h^β / Gamma[β+2] *Table[Which[i<j, 0, i==j, 1, i>j,
    (i-j+1)^(β+1) - 2(i-j)^(β+1) + 2(i-j-1)^(β+1)], {i,m}, {i,m}];
P05:=H[0.5,h,m]; P07:=H[0.7,h,m]; P1:=H[1,h,m]; P2:=H[2,h,m];
```

- *Формирование правой части уравнения в операционной области:*

```
Φ:=F-y1*D1.P07.One-D2.P2.(y0*One+y1+Tim);
```

- *Нахождение изображения решения в операционной области (функции, первой и второй производной):*

```
U:=Inverse[(P05+D1.P07+D2.P2)].Φ;
Y1:=y1*One+P1.U;
Y:=y0*One+y1*Tim+P2.U;
```

- *Визуализация векторов решения в операционной области (рис. 1):*

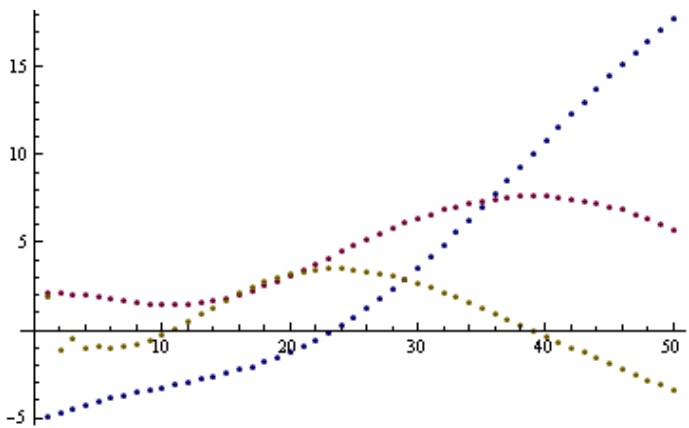

***Рис. 1.*** *Элементы векторов* **U, Y1, Y**

- *Формирование аппроксимаций решения уравнения и его двух производных:*
  ```
  ya=Y.V; y1a=Y1.V;
  ```
- *Визуализация аппроксимаций решения уравнения и его двух производных (рис.2):*
  ```
  Plot[{ya,y1a,ua}, {t,0,5}, Plot Points ->200]
  ```

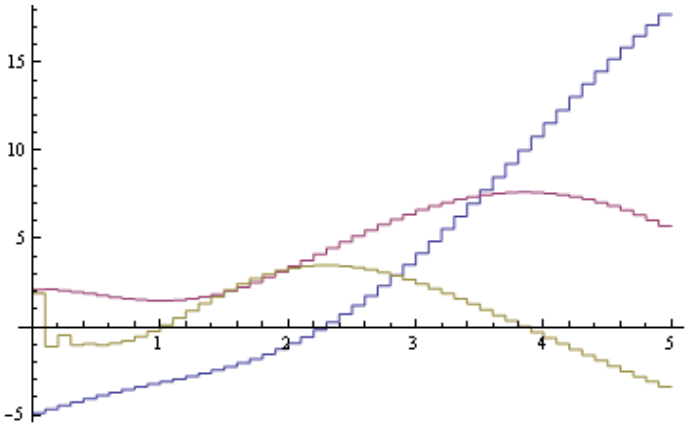

***Рис.2.*** *Аппроксимация решения уравнения (41).*

**Заключение.** Предложен порядок применения операционного метода на основе полиномиальной аппроксимации сигналов к решению линейных интегро-дифференциальных уравнений дробного и смешанного порядков с производными по Капуто. Порядок предполагает замену производной наивысшего целого порядка новой функцией, для которой исходное дифференциальное уравнение превращается в интегральное с операторами Римана–Лиувилля. В операционной области используется блочно-импульсная система базисных функций, для которой существенно упрощается операция произведения функций. Рассмотренный метод может быть реализован в любой из известных систем компьютерной алгебры.

### Список использованной литературы:

An application of approximated operational approach to the solution of linear non-integer order differential equations with variable coefficients has been considered. It has been shown, that the method can be applied also to the solution both of the non-integer and combined order differential equations. The computer experiments were fulfilled in Mathematica® program area.

**Key words:** *Polynomial approximation, fractional calculus, non-integer order differential equation, Caputo's derivative, mathematical modelling, dynamical system.*